\documentclass[11pt]{article}

\usepackage{tikz}
\usepackage[margin=1in]{geometry}
\usepackage{proof}
\usepackage{amsmath,amsthm,amssymb}
\usepackage{csquotes}
\usepackage{amsthm}
\usepackage{pifont}
\usepackage{graphicx}
\usepackage{setspace}
\usepackage{enumitem}
\usepackage{epigraph}
\usepackage{multicol}
\usepackage{xfrac}
\usepackage{cancel}
\usepackage{wasysym}
\usepackage{pxfonts}
\usepackage{textcomp}
\usepackage{textpos}
\usepackage{tikz-cd}
\usepackage{array}
\usepackage[authoryear]{natbib}
\usepackage{booktabs}
\usepackage{tabularx}
\usepackage{float}
\usepackage[hidelinks]{hyperref}

\newcommand{\assert}[1]{+\langle #1\rangle}
\newcommand{\deny}[1]{-\langle #1\rangle}
\newcommand{\stroke}{\mathbin{\mid}}

\newcommand{\signa}{\textit{\textbf{a}}}
\newcommand{\signb}{\textit{\textbf{b}}}
\newcommand{\signc}{\textit{\textbf{c}}}
\newcommand{\signd}{\textit{\textbf{d}}}
\newcommand{\signe}{\textit{\textbf{e}}}
\newcommand{\signf}{\textit{\textbf{f}}}
\newcommand{\signg}{\textit{\textbf{g}}}

\newcommand{\nc}{\,\mid\!\sim\,}

\newcommand{\rcoimp}{\mathbin{\rotatebox[origin=c]{180}{$\coimp$}}}

\newcommand{\coimp}{>\hspace{-.2cm}-\hspace{-.2cm}-}

\newcommand{\starredconnective}[2]{%
  \mathbin{\vbox{\offinterlineskip
    \halign{\hfil$##$\hfil\cr
      \scriptstyle *\cr\noalign{\kern#1}#2\cr}}}%
}
\newcommand{\starimp}{%
  \starredconnective{-1.3pt}{\rightarrow}%
}

\newcommand{\starcoimp}{%
  \starredconnective{-1.5pt}{\rcoimp}%
}

\newcommand{\coimprule}{\mbox{\(\coimp\)}}

\newcommand{\csim}{\sim \hspace{-.1cm} }

\makeatletter
\newcommand*\bigcdot{\mathpalette\bigcdot@{.8}}
\newcommand*\bigcdot@[2]{\mathbin{\vcenter{\hbox{\scalebox{#2}{$\m@th#1\bullet$}}}}}
\makeatother



\makeatletter
\providecommand*{\dashv}{
  \mathrel{
    \mathpalette\@dashv\nc
  }
}
\newcommand*{\@dashv}[2]{
  \reflectbox{$\m@th#1#2$}
}
\makeatother

\title{On the Definability of Strong Negation\\  in Bilateral Logics}
\author{Ryan Simonelli\\ \href{https://orcid.org/0000-0003-3754-1671}{ORCID: 0000-0003-3754-1671} \\ \\ Zhejiang University, School of Philosophy}

\begin{document}

\maketitle

\begin{textblock*}{145mm}(.01\textwidth,-7.8cm)
	\small \textbf{Penultimate draft.  Please cite the published version in \textit{The Review of Symbolic Logic}, once that version is available.}
\end{textblock*}

\begin{abstract}
Recent work in bilateral logic has concerned itself with whether and to what extent strong negation, which ``toggles'' between assertion and denial, is definable in bilateral systems that lack it.  This paper presents a number of results about the definability of strong negation in bilateral systems.  First, I show that a constructive Nelson-style Sheffer stroke defines the strong negation of $A$ as $(A\stroke(A\stroke A)) \stroke A$.  The same stroke defines the constructive Nelson implication, thus providing a single-connective basis for the \(\{\mathord{\sim},\mathord{\to}\}\)-fragment of N4. Second, I provide an exhaustive characterization of the eight combinations of a single ``aggregative'' connective and a single Nelson-style connective that suffice to jointly define strong negation, with neither defining it individually.  One such combination is the material conditional along with constructive co-implication, with which the strong negation of $A$ is definable simply as $A \supset (A \rcoimp A)$.  Third, I show that a constructive and connexive Wansing-style Sheffer stroke (recently referred to as ``connexive exclusion'') likewise defines the strong negation of $A$ as $(A\stroke(A\stroke A)) \stroke A$ and also defines connexive implication.  Fourth, I show that, unlike with the Nelson-style connectives, there are no combinations of Wansing-style and aggregative connectives that jointly define strong negation if neither connective defines it individually.\\

\noindent \textbf{Keywords:} Bilateralism, Strong Negation, Nelson's N4, Connexive Logic

\noindent\textbf{Mathematics Subject Classification:} 03B60; 03B53, 03F03.
\end{abstract}

\section{Introduction}

Bilateral proof systems provide both positive and negative rules for connectives, typically interpreted as rules for proof and refutation, or, alternatively, assertion and denial.   Since their initial development three decades ago in the work of \cite{Smiley1996} and \cite{Rumfitt2000}, the defining feature of such systems has generally been taken to be a strong negation operator that ``toggles'' between proof and refutation, or, alternatively, between assertion and denial.  The possibility of such negation rules in bilateral systems makes them a natural setting for constructive logics with strong negation in the tradition of Nelson \citeyearpar{Nelson1949}, such as his paraconsistent N4 \citep{AlmukdadNelson1984}.\footnote{N4 is the paraconsistent variant of the logic introduced in \citep{Nelson1949}. See \cite{KamideWansing2012} for a proof-theoretic overview of such logics.}  Indeed, as \cite{Gibbard2002} points out, Rumfitt's bilateral natural deduction system for classical logic, if we leave out the ``coordination principles'' which impose the exclusivity and exhaustivity of assertion and denial,  in fact just is a bilateral system for N4.  Though toggle negation is traditionally the hallmark of bilateral systems, in recent years, there has been substantial interest in bilateral systems that lack a primitive toggle negation operator.\footnote{See, for instance, \cite{Francez2014}, \cite{Drobyshevich2021},\cite{WansingNikiDrobyshevich2025}, \cite{Ayhan2025}, and \cite{Ayhan2026}.  See also \cite{Ayhan2025b} for some philosophical motivation, from a feminist perspective (drawing on \cite{Plumwood1993}), for not wanting toggle negation to be definable in a bilateral system.  }  In the context of such systems, a natural question is whether toggle negation, though not primitive, is nevertheless definable.  Two recent results, one positive and one negative, are particularly worth noting here.

Positively, \cite{Oddsson2025} has shown that strong negation is in fact definable in the notable bilateral intuitionistic logic \(2\mathrm{Int}\), introduced by \cite{Wansing2016}.  Where $\rightarrow$ is constructive implication and $\rcoimp$ is the constructive co-implication operator belonging to 2Int (to be explained below), Oddsson shows that toggle negation is definable as follows:

\[
 \mathord{\sim} A :=
 ((A\wedge(A\to(A \rcoimp A))))
 \vee ((A\to A) \rcoimp A).
\]

\noindent Thus, insofar as a bilateral system for Nelson's N4 can be obtained by adding toggle negation to 2Int, Oddsson's result shows that 2Int by itself already implicitly contains Nelson's N4. Negatively, Wansing, Niki, and Drobyshevich \citeyearpar{WansingNikiDrobyshevich2025} have shown that toggle negation is  \textit{not} definable in the logic that is reasonably regarded as the \textit{connexive} variant of 2Int. Their bi-connexive logic B2C contains the same four binary connectives as 2Int---the ``aggregative'' connectives of conjunction and disjunction, and the constructive conditionals of implication and co-implication.  The two conditionals, however, are made connexive through a variation of their deniability and assertibility conditions, respectively.  Unlike 2Int, strong negation is not definable in B2C.

These recent results raise a number of questions regarding the definability of strong negation in bilateral constructive logics.  Oddsson's definition of toggle negation in 2Int uses all four of the binary connectives typically included in 2Int, and it is natural to wonder whether a negation-free constructive logic with a smaller binary connective toolkit can still define toggle negation, and, if so, which ones do. The first pair of results here answers these questions.  First, I show that a system consisting in a single constructive Sheffer stroke suffices to define strong negation along with constructive implication, and thus, suffices to define the $\{\sim, \rightarrow\}$ fragment of N4.   Second, I show that there are a number of negation-free systems containing two connectives, neither of which is sufficient to define strong negation on their own, which jointly suffice for the entirety of N4.  Among them is the system containing the material conditional and constructive co-implication.  Wansing, Niki, and Drobyshevich's result, on the other hand, shows that negation is not definable using the standard connectives of constructive connexive logic, and it is natural to wonder whether constructive connexive connectives can be used to define toggle negation at all, either by themselves or with aggregative connectives.  The second pair of results answers these questions.  First, I show that a system consisting in a single constructive and connexive Sheffer stroke, recently discussed by \cite{ShramkoWansing2025} under the name of ``connexive exclusion,'' suffices to define strong negation along with connexive implication.    Second, I show that unlike with the Nelson-style connectives, there are no combinations of Wansing-style and aggregative connectives that jointly define strong negation if neither connective defines it individually.

\section{A Generalized Approach to Bilateral Proof Systems}

Let $\mathcal{A}$ be a stock of atomic formulas, and let $\mathcal{L}$ be the language generated from $\mathcal{A}$ by the connectives under consideration. A \textit{signed formula} is an expression of the form
$+\langle A\rangle$ or $-\langle A\rangle$, where $A\in\mathcal{L}$.  $+\langle A\rangle$ expresses the assertion of $A$ and $-\langle A\rangle$ expresses the denial of $A$.  These signs belong to the proof-theoretic metalanguage, and are to be sharply distinguished from object-language operators such as nullation or negation.  In particular, unlike a negation operator, they can be neither embedded in sentences nor iterated. Let $\mathcal{L}_{\pm} = \{+\langle A\rangle\mid A\in\mathcal{L}\} \cup \{-\langle A\rangle\mid A\in\mathcal{L}\}.$  The metavariables $\varphi$ and $\psi$ range over signed formulas, while $\Gamma$ and $\Delta$ range over finite multisets of signed formulas.   I will work in a multiple conclusion bilateral sequent calculus.\footnote{Nothing hangs on this choice, and all of the results that follow may equally be obtained in a bilateral natural deduction setting, using the generalized natural deduction rules I propose in \citep{Simonelli2025}, but omitting the coordination principles.  In particular, variants on all of the results below are easily obtained using the BNK1 rules for the aggregative connectives,  BNK2 rules for the Nelson-style connectives, and the natural connexive variants of the BNK2 rules for the Wansing-style connectives.}  As I've proposed elsewhere \citeyearpar{Simonelli2026b}, a sequent of the form $\Gamma \vdash \Delta$ might be read in a way inspired by \cite{Restall2005} as saying that \textit{making} all of the assertions and denials in $\Gamma$ along with \textit{challenging} all of the assertions and denials in $\Delta$ is incoherent or ``out of bounds.''

Instead of defining connectives one by one, I will define them by way of general rule schemas, following the approach I've developed elsewhere \citeyearpar{Simonelli2024,Simonelli2025,Simonelli2026}. Most recently, in proposing bilateral systems for the logics in the FDE family, I \citeyearpar{Simonelli2026} consider binary connectives of the following form:\footnote{Rules of this form are standard in sequent calculi for Belnap-Dunn logics, though the negative rules are typically presented with a negation operator rather than a denial sign.  See, for instance, \cite{ArieliAvron1998}.  \cite{Wintein2016} presents 4-signed rules of this form for logics in the FDE family, which are inter-translatable with the bilateral rules presented here.}

\begin{quote}
\textbf{Aggregative Connectives:}

\begin{center}
   \begin{multicols}{2}
  
       $
  \infer[{{\textbf{\textit{c}}}_{\circ_\text{L}}}]{\Gamma, \textbf{\textit{c}} \langle A \circ B \rangle \vdash \Delta}{\Gamma, \textbf{\textit{a}}  \langle A \rangle, \textbf{\textit{b}} \langle B \rangle \vdash \Delta}
  $ 
  
  \columnbreak
  
         $
  \infer[{{\textbf{\textit{c}}}_{\circ_\text{R}}}]{\Gamma \vdash  \textbf{\textit{c}}  \langle A \circ B \rangle, \Delta}{\Gamma \vdash  \textbf{\textit{a}}  \langle  A \rangle, \Delta & \Gamma \vdash \textbf{\textit{b}} \langle B \rangle, \Delta}
  $
    
   \end{multicols}

   \begin{multicols}{2}
   
    $
  \infer[{{\textbf{\textit{c}}^*}_{\circ_\text{L}}}]{\Gamma,  \textbf{\textit{c}}^*  \langle A \circ B \rangle \vdash \Delta}{\Gamma, \textbf{\textit{a}}^*  \langle  A \rangle \vdash \Delta & \Gamma, \textbf{\textit{b}}^* \langle B \rangle \vdash \Delta}
  $
  
    \columnbreak
  
       $
  \infer[{{\textbf{\textit{c}}^*}_{\circ_\text{R}}}]{\Gamma \vdash \textbf{\textit{c}}^* \langle A \circ B \rangle, \Delta}{\Gamma \vdash \textbf{\textit{a}}^*  \langle A \rangle,  \textbf{\textit{b}}^* \langle B \rangle, \Delta}
  $

   \end{multicols}
\end{center}
\end{quote}

\noindent Here, \textbf{\textit{a}}, \textbf{\textit{b}}, and \textbf{\textit{c}} are variables ranging over $\{+, -\}$ and $^*$ is a function mapping $+$ to $-$ and $-$ to $+$.  Using this schema, we can define eight aggregative connectives, considering every possible assignment of signs to \textbf{\textit{a}}, \textbf{\textit{b}}, and \textbf{\textit{c}}.  In particular, we can define the following connectives:

\begin{center}
\begin{multicols}{2}

$\wedge_\circ$: $\textbf{\textit{a}}=+, \textbf{\textit{b}}=+, \textbf{\textit{c}}=+$\\

$\mid_\circ$: $\textbf{\textit{a}}=+, \textbf{\textit{b}}=+, \textbf{\textit{c}}=- $\\

$\rightarrow_\circ$: $\textbf{\textit{a}}=+, \textbf{\textit{b}}=-, \textbf{\textit{c}}=-$\\

$\starimp_\circ$: $\textbf{\textit{a}}=+, \textbf{\textit{b}}=-, \textbf{\textit{c}}=+$

\columnbreak

$\vee_\circ$: $\textbf{\textit{a}}=-, \textbf{\textit{b}}=-, \textbf{\textit{c}}=-$\\

$\downarrow_\circ$: $\textbf{\textit{a}}=-, \textbf{\textit{b}}=-, \textbf{\textit{c}}=+$\\

$\coimp_\circ$: $\textbf{\textit{a}}=-, \textbf{\textit{b}}=+, \textbf{\textit{c}}=+$\\

$\starcoimp_\circ$:$\textbf{\textit{a}}=-, \textbf{\textit{b}}=+, \textbf{\textit{c}}=-$

\end{multicols}
\end{center}

\noindent  Of course, we have here the common connectives: conjunction, disjunction, and the material conditional, as well as some less common but still familiar ones---the Sheffer stroke and Peirce's arrow. We also have the ``anti-conditional,'' $\starimp$, equivalent to the conditional's strong negation, as well as the De Morgan dual of the material conditional, $\coimp$, and its opposite, $\starcoimp$.   Note that these are understood as potentially primitive binary operators, and are definable in a bilateral system without any use of a unary negation operator.  Hence, I write $\starimp$ rather than $\nrightarrow$ to be clear that this connective is not to be understood as the negation of $\rightarrow$, even though it is inferentially equivalent to it.

Now, in a bilateral sequent calculus of this sort, ``toggle negation'' will be given the obvious rules, toggling between assertion and denial:

\begin{center}

\begin{multicols}{4}

$\infer[+_{\sim_L}]{\Gamma, + \langle \csim A \rangle \vdash \Delta}{\Gamma, - \langle A \rangle \vdash \Delta}$

\columnbreak

$\infer[-_{\sim_L}]{\Gamma, - \langle \csim A \rangle \vdash \Delta}{\Gamma, + \langle A \rangle \vdash \Delta}$

\columnbreak

$\infer[+_{\sim_R}]{\Gamma \vdash  + \langle \csim A \rangle, \Delta}{\Gamma \vdash  - \langle A \rangle, \Delta}$

\columnbreak

$\infer[-_{\sim_R}]{\Gamma \vdash  - \langle \csim A \rangle, \Delta}{\Gamma \vdash  + \langle A \rangle, \Delta}$

\end{multicols}
\end{center}

\noindent However, given some of the binary connectives defined above, we need not take toggle negation as primitive.  Consider, for instance, the system consisting in just the rules for the Sheffer stroke (to avoid clutter I omit the subscripted $\circ$ here):

\begin{quote}
\begin{center}

\begin{multicols}{2}

$
\infer[+_{\stroke\mathrm L}]
      {\Gamma,\assert{A\stroke B}\vdash\Delta}
      {\Gamma,\deny A\vdash\Delta
       &
       \Gamma,\deny B\vdash\Delta}
$

\columnbreak

$
\infer[+_{\stroke\mathrm R}]
      {\Gamma\vdash\assert{A\stroke B},\Delta}
      {\Gamma\vdash\deny A,\deny B,\Delta}
$

\end{multicols}

\begin{multicols}{2}

$
\infer[-_{\stroke\mathrm L}]
      {\Gamma,\deny{A\stroke B}\vdash\Delta}
      {\Gamma,\assert A,\assert B\vdash\Delta}
$

\columnbreak

$
\infer[-_{\stroke\mathrm R}]
      {\Gamma\vdash\deny{A\stroke B},\Delta}
      {\Gamma\vdash\assert A,\Delta
       &
       \Gamma\vdash\assert B,\Delta}
$

\end{multicols}

\end{center}
\end{quote}

\noindent Using this connective, we can define toggle negation in the usual way as $A \mid A$.  Thus, given the axiom schema of Contexted Reflexivity and the Rule of Cut:

\begin{center}
\begin{multicols}{2}
$\infer[_\text{Contexted Reflex.}]{\Gamma, \varphi \vdash \varphi, \Delta}{}$

\columnbreak

$\infer[_\text{Cut}]{\Gamma, \Gamma' \vdash \Delta, \Delta'}{\Gamma, \varphi \vdash \Delta & \Gamma' \vdash \varphi, \Delta'}$
\end{multicols}
\end{center}

\noindent We can derive all of the rules for toggle negation as follows:

\begin{center}
\begin{multicols}{2}
\small

$
\infer[_\mathrm{Cut}]
{
  \Gamma,\assert{A\stroke A}
  \vdash
  \Delta
}
{
  \infer[+_{\stroke\mathrm L}]
  {
    \assert{A\stroke A}
    \vdash
    \deny A
  }
  {
    \infer[_\text{C. Reflex.}]
    {
      \deny A
      \vdash
      \deny A
    }
    {}
    &
    \infer[_\text{C. Reflex.}]
    {
      \deny A
      \vdash
      \deny A
    }
    {}
  }
  &
  \Gamma,\deny A
  \vdash
  \Delta
}
$

\columnbreak 

$
\infer[_\mathrm{Cut}]
{
  \Gamma
  \vdash
  \assert{A\stroke A},\Delta
}
{
  \Gamma
  \vdash
  \deny A,\Delta
  &
  \infer[+_{\stroke\mathrm R}]
  {
    \deny A
    \vdash
    \assert{A\stroke A}
  }
  {
    \infer[_\text{C. Reflex.}]
    {
      \deny A
      \vdash
      \deny A,\deny A
    }
    {}
  }
}
$
\end{multicols}
\begin{multicols}{2}
\small
$
\infer[_\mathrm{Cut}]
{
  \Gamma,\deny{A\stroke A}
  \vdash
  \Delta
}
{
  \infer[-_{\stroke\mathrm L}]
  {
    \deny{A\stroke A}
    \vdash
    \assert A
  }
  {
    \infer[_\text{C. Reflex.}]
    {
      \assert A,\assert A
      \vdash
      \assert A
    }
    {}
  }
  &
  \Gamma,\assert A
  \vdash
  \Delta
}
$

\columnbreak

$
\infer[_\mathrm{Cut}]
{
  \Gamma
  \vdash
  \deny{A\stroke A},\Delta
}
{
  \Gamma
  \vdash
  \assert A,\Delta
  &
  \infer[-_{\stroke\mathrm R}]
  {
    \assert A
    \vdash
    \deny{A\stroke A}
  }
  {
    \infer[_\text{C. Reflex.}]
    {
      \assert A
      \vdash
      \assert A
    }
    {}
    &
    \infer[_\text{C. Reflex.}]
    {
      \assert A
      \vdash
      \assert A
    }
    {}
  }
}
$
\end{multicols}
\end{center}

\noindent  Indeed, given just these rules, we have a single-connective bilateral system for FDE.\footnote{I show this in \citep{Simonelli2026}. I'm not aware of a prior formulation of a single connective system for FDE.}  Likewise for the single connective system containing the rules for Peirce's arrow.

I have referred to connectives of the above form as \textit{aggregative connectives}, since they simply function to collect formulas of the different signs on the sides of the turnstile to which those formulas belong.\footnote{The terminology is from Brandom (see \cite{HlobilBrandom2024}).  Brandom refers to conjunction and disjunction as ``Boolean helper-monkeys.'' A better term for all eight aggregative connectives given by this bilateral schema might be ``Dunnean helper-monkeys.''}  They do not internalize inferential relations across the turnstile.  Thus, though we have a conditional operator, $\rightarrow_\circ$, it is the material conditional of FDE, and, without the addition of further coordination principles imposing the exclusivity or exhaustivity of assertion and denial, it does not validate modus ponens $+ \langle A \rangle, + \langle A \rightarrow_\circ B \rangle \vdash + \langle B \rangle$ or identity $\vdash + \langle A \rightarrow_\circ A \rangle$.  It is thus natural to want to add a conditional that does.  Here, I will consider the \textit{constructive Nelson conditional}, which figures in Nelson's N4 as well as Wansing's 2Int.\footnote{I focus on the constructive variants of these conditionals here, since the recent literature on the definability of strong negation in bilateral systems has mainly dealt with constructive bilateral systems.  The main results here carry over to non-constructive variants of these conditionals, rules for which are obtained simply by removing the restriction to single conclusions in the rules and eliminating the repetition of the conditional in the premise of the left rule.} Though the Nelson conditional is typically presented in a natural deduction system or single conclusion sequent calculus,  \cite{Drobyshevich2026} gives the following rules for the constructive Nelson conditional in a multiple conclusion bilateral sequent calculus:

 \begin{center}

   \begin{multicols}{2}
   
       $
  \infer[+_{\rightarrow_\text{L}}]
  {\Gamma, + \langle A \rightarrow B \rangle \vdash \Delta}
  {\Gamma, + \langle A \rightarrow B \rangle
  \vdash + \langle A \rangle
  &
  \Gamma, + \langle B \rangle \vdash \Delta}
  $
  
   \columnbreak
  
       $
  \infer[+_{\rightarrow_\text{R}}]
  {\Gamma \vdash + \langle A \rightarrow B \rangle, \Delta}
  {\Gamma, + \langle A \rangle
  \vdash + \langle B \rangle}
  $
  
   \end{multicols}

   \begin{multicols}{2}
  
       $
  \infer[-_{\rightarrow_\text{L}}]
  {\Gamma, - \langle A \rightarrow B \rangle \vdash \Delta}
  {\Gamma, + \langle A \rangle, - \langle B \rangle \vdash \Delta}
  $ 
  
  \columnbreak
  
       $
  \infer[-_{\rightarrow_\text{R}}]
  {\Gamma \vdash - \langle A \rightarrow B \rangle, \Delta}
  {\Gamma \vdash + \langle A \rangle, \Delta
  &
  \Gamma \vdash - \langle B \rangle, \Delta}
  $
    
   \end{multicols}

\end{center}

\noindent Note the restriction to single conclusions in the premise of the positive right rule and in the left premise of the positive left rule.  This restriction ensures that the conditional behaves constructively (the duplicated formula in the positive left rule ensures that Contraction is admissible).  

Now, continuing with the same generalized approach to bilateralism, we need not restrict ourselves to this one conditional, but, rather, we can consider the whole range of connectives with rules of this general form:

\begin{quote}
\textbf{Constructive Nelson-Style Connectives:}

\begin{center}

\begin{multicols}{2}

$
\infer[{{\textbf{\textit{c}}}_{\triangleright_\mathrm{L}}}]
{\Gamma,
 \textbf{\textit{c}}\langle A\triangleright B\rangle
 \vdash\Delta}
{\Gamma,
 \textbf{\textit{a}}\langle A\rangle,
 \textbf{\textit{b}}\langle B\rangle
 \vdash\Delta}
$

\columnbreak

$
\infer[{{\textbf{\textit{c}}}_{\triangleright_\mathrm{R}}}]
{\Gamma\vdash
 \textbf{\textit{c}}\langle A\triangleright B\rangle,
 \Delta}
{\Gamma\vdash
 \textbf{\textit{a}}\langle A\rangle,\Delta
 &
 \Gamma\vdash
 \textbf{\textit{b}}\langle B\rangle,\Delta}
$

\end{multicols}

\begin{multicols}{2}

$
\infer[{{\textbf{\textit{c}}^*}_{\triangleright_\mathrm{L}}}]
{\Gamma,
 \textbf{\textit{c}}^*
 \langle A\triangleright B\rangle
 \vdash\Delta}
{\Gamma,
 \textbf{\textit{c}}^*
 \langle A\triangleright B\rangle
 \vdash
 \textbf{\textit{a}}\langle A\rangle
 &
 \Gamma,
 \textbf{\textit{b}}^*\langle B\rangle
 \vdash\Delta}
$

\columnbreak

$
\infer[{{\textbf{\textit{c}}^*}_{\triangleright_\mathrm{R}}}]
{\Gamma\vdash
 \textbf{\textit{c}}^*
 \langle A\triangleright B\rangle,
 \Delta}
{\Gamma,
 \textbf{\textit{a}}\langle A\rangle
 \vdash
 \textbf{\textit{b}}^*\langle B\rangle}
$

\end{multicols}

\end{center}
\end{quote}

\noindent As before, we can define eight distinct connectives, each determined by some assignment of signs to  $\textbf{\textit{a}}$, $\textbf{\textit{b}}$, and $\textbf{\textit{c}}$:

\begin{center}
\begin{multicols}{2}

$\wedge_\triangleright$: $\textbf{\textit{a}}=+, \textbf{\textit{b}}=+, \textbf{\textit{c}}=+$\\

$\mid_\triangleright$: $\textbf{\textit{a}}=+, \textbf{\textit{b}}=+, \textbf{\textit{c}}=- $\\

$\rightarrow_\triangleright$: $\textbf{\textit{a}}=+, \textbf{\textit{b}}=-, \textbf{\textit{c}}=-$\\

$\starimp_\triangleright$: $\textbf{\textit{a}}=+, \textbf{\textit{b}}=-, \textbf{\textit{c}}=+$

\columnbreak

$\vee_\triangleright$: $\textbf{\textit{a}}=-, \textbf{\textit{b}}=-, \textbf{\textit{c}}=-$\\

$\downarrow_\triangleright$: $\textbf{\textit{a}}=-, \textbf{\textit{b}}=-, \textbf{\textit{c}}=+$\\

$\coimp_\triangleright$: $\textbf{\textit{a}}=-, \textbf{\textit{b}}=+, \textbf{\textit{c}}=+$\\

$\starcoimp_\triangleright$:$\textbf{\textit{a}}=-, \textbf{\textit{b}}=+, \textbf{\textit{c}}=-$

\end{multicols}
\end{center}

\noindent Attending to the subscripts, we can distinguish between the material conditional $\rightarrow_\circ$ and the constructive Nelson conditional, $\rightarrow_\triangleright$.  Likewise, we obtain the familiar constructive co-implication ${\coimp}_\triangleright$, modulo a harmless reversal of argument order: what is standardly written $A \rcoimp B$ will here be written $B \coimp A$. These two constructive implication operators are standard in constructive bilateral logics.  We can also consider constructive conjunction and disjunction connectives, $\wedge_\triangleright$ and $\vee_\triangleright$, which might be understood as the Nelson-style analogues of fusion and fission (intensional conjunction and disjunction), familiar from relevance logic, obtained by instantiating their familiar definitions with the constructive Nelson conditional.\footnote{In fact, the bilateral natural deduction system for classical logic proposed by \cite{delValle-Inclan2023}, and which I discuss under the name ``BNK2'' \citep{Simonelli2025}, can be understood as having rules for these operators, though, like the conditional in Rumfitt's natural deduction system, their constructivity (and, indeed, asymmetry) is destroyed in del Valle-Inclan and Schloeder's system by the inclusion of coordination principles.}  The first object of interest here, however, will be the  \emph{constructive
Nelson-style Sheffer stroke}.

\section{Defining Strong Negation with One Connective}

The Nelson-style Sheffer stroke is given by the following rules:\footnote{Though this connective is extremely natural in a bilateral context, as far as I'm aware, this is the first time it has been explicitly considered.  One seemingly related construction is the ``Sheffer stroke Nelson algebras'' considered in  \cite{OnerKaticanSaeid2025}.  However, with the Sheffer identities assumed throughout their proofs, the algebras they define are in fact Boolean! In particular, their Lemma~3.2 establishes a bounded distributive lattice, and the verification of (N6) in Theorem~3.1 establishes the Boolean complement laws. Their construction therefore yields the classical Sheffer stroke.}

\begin{center}

\begin{multicols}{2}

$\infer[+_{\stroke\mathrm L}]
      {\Gamma,\assert{A\stroke B}\vdash\Delta}
      {\Gamma,\assert{A\stroke B}\vdash\assert A
       &
       \Gamma,\deny B\vdash\Delta}$

\columnbreak

$\infer[+_{\stroke\mathrm R}]
      {\Gamma\vdash\assert{A\stroke B},\Delta}
      {\Gamma,\assert A\vdash\deny B}$

\end{multicols} 

\begin{multicols}{2}

$\infer[-_{\stroke\mathrm L}]
      {\Gamma,\deny{A\stroke B}\vdash\Delta}
      {\Gamma,\assert A,\assert B\vdash\Delta}$

\columnbreak

$\infer[-_{\stroke\mathrm R}]
      {\Gamma\vdash\deny{A\stroke B},\Delta}
      {\Gamma\vdash\assert A,\Delta
       &
       \Gamma\vdash\assert B,\Delta}$
\end{multicols}

\end{center}

\noindent Thus, denying $A \stroke B$ has exactly the same conditions as asserting an ordinary conjunction: it amounts to asserting both $A$ and $B$. The condition for \textit{asserting} $A \stroke B$, however, is that the denial of $B$ follow from the assertion of $A$.  In the context of a bilateral system, this is a perfectly reasonable connective to include as a primitive binary connective, indeed, as one's only primitive connective.  We will now show the following:

\begin{quote}
\textbf{Proposition 1:} The system containing just the Nelson-style Sheffer stroke rules enables one to define strong negation and thus, all Nelson-style connectives. 
\end{quote}

\textit{Proof:} We define strong negation as follows: 

\begin{quote}
$\csim A := (A \mid (A \mid A)) \mid A$
\end{quote}

\noindent To show that this indeed defines $\csim A$, we show that $+\langle (A \mid (A \mid A)) \mid A \rangle \dashv\vdash -\langle A\rangle$ and $-\langle (A \mid (A \mid A)) \mid A \rangle \dashv\vdash +\langle A\rangle$:

\begin{quote}
\small 

$-\langle A\rangle \vdash +\langle (A \mid (A \mid A)) \mid A \rangle$:

\begin{center}
$\infer[+_{\mid_{\mathrm R}}]
{
  -\langle A\rangle
  \vdash
  +\langle (A \mid (A \mid A)) \mid A\rangle
}
{
  \infer[_\text{C. Reflex.}]
  {
    -\langle A\rangle,
    +\langle A \mid (A \mid A)\rangle
    \vdash
    -\langle A\rangle
  }
  {}
}$
\end{center}

$+\langle (A \mid (A \mid A)) \mid A \rangle \vdash -\langle A\rangle $:

\begin{center}

$\infer[+_{\mid_{\mathrm L}}]
{
  +\langle (A \mid (A \mid A)) \mid A\rangle
  \vdash
  -\langle A\rangle
}
{
  \infer[+_{\mid_{\mathrm R}}]
  {
    +\langle (A \mid (A \mid A)) \mid A\rangle
    \vdash
    +\langle A \mid (A \mid A)\rangle
  }
  {
    \infer[-_{\mid_{\mathrm R}}]
    {
      +\langle (A \mid (A \mid A)) \mid A\rangle,
      +\langle A\rangle
      \vdash
      -\langle A \mid A\rangle
    }
    {
      \infer[_\text{C. Reflex.}]
      {
        +\langle (A \mid (A \mid A)) \mid A\rangle,
        +\langle A\rangle
        \vdash
        +\langle A\rangle
      }
      {}
      &
      \infer[_\text{C. Reflex.}]
      {
        +\langle (A \mid (A \mid A)) \mid A\rangle,
        +\langle A\rangle
        \vdash
        +\langle A\rangle
      }
      {}
    }
  }
  &
  \infer[_\text{C. Reflex.}]
  {
    -\langle A\rangle
    \vdash
    -\langle A\rangle
  }
  {}
}$
\end{center}

$+\langle A\rangle \vdash -\langle (A \mid (A \mid A)) \mid A \rangle$:

\begin{center}
$\infer[-_{\mid_{\mathrm R}}]
{
  +\langle A\rangle
  \vdash
  -\langle (A \mid (A \mid A)) \mid A\rangle
}
{
  \infer[+_{\mid_{\mathrm R}}]
  {
    +\langle A\rangle
    \vdash
    +\langle A \mid (A \mid A)\rangle
  }
  {
    \infer[-_{\mid_{\mathrm R}}]
    {
      +\langle A\rangle,
      +\langle A\rangle
      \vdash
      -\langle A \mid A\rangle
    }
    {
      \infer[_\text{C. Reflex.}]
      {
        +\langle A\rangle,
        +\langle A\rangle
        \vdash
        +\langle A\rangle
      }
      {}
      &
      \infer[_\text{C. Reflex.}]
      {
        +\langle A\rangle,
        +\langle A\rangle
        \vdash
        +\langle A\rangle
      }
      {}
    }
  }
  &
  \infer[_\text{C. Reflex.}]
  {
    +\langle A\rangle
    \vdash
    +\langle A\rangle
  }
  {}
}$

\end{center}

$-\langle (A \mid (A \mid A)) \mid A \rangle \vdash +\langle A\rangle $:

\begin{center}

$
\infer[-_{\mid_{\mathrm L}}]
{
  -\langle (A \mid (A \mid A)) \mid A\rangle
  \vdash
  +\langle A\rangle
}
{
  \infer[_\text{C. Reflex.}]
  {
    +\langle A \mid (A \mid A)\rangle,
    +\langle A\rangle
    \vdash
    +\langle A\rangle
  }
  {}
}
$
\end{center}
\normalsize
\end{quote}

\noindent Given Cut, we can derive the standard positive and negative left and right rules for toggle negation.  

 Having defined strong negation, we may now define the constructive Nelson conditional as follows:\footnote{Thus, using only the Sheffer stroke, the conditional $A \rightarrow B$ may be defined as $A\stroke(((B\stroke(B\stroke B))\stroke B))$.  A much shorter stroke-only definition of $A \rightarrow B$ as $A\stroke(A\stroke B)$ is possible, though I omit the proofs here.}
\begin{quote}
$A\rightarrow B := A\stroke\csim B$
\end{quote}

\noindent To see why this definition works, substitute $\csim B$ for $B$ in the rules for $\stroke$. The derived rules for strong negation make $\assert{\csim B}$ interchangeable with $\deny B$, and $\deny{\csim B}$ interchangeable with $\assert B$. Consequently, the positive rules for $A\stroke\csim B$ are exactly the positive rules for the constructive Nelson conditional, while its negative rules are exactly the negative conditional rules. Thus, the Nelson-style Sheffer stroke defines both strong negation and constructive implication. The remaining Nelson-style connectives may be similarly defined. It follows that $\stroke$ is by itself expressively complete for the family of Nelson-style connectives considered here.  $\Box$

In order to move on to the next main result, we will need to introduce some basic semantic machinery.  

\section{A Generalized Approach to Bilateral Semantics}

Let us start with the notion of validity.  Once again, I read $\Gamma \vdash \Delta$ as saying that making all of the moves in $\Gamma$ and challenging all of the moves in $\Delta$ is incoherent or ``out of bounds.''  That is to say, there is no way for all of the moves in $\Gamma$ to be correct and all of the moves in $\Delta$ to be incorrect.  Intuitively, we might regard an assertion as correct just in case the sentence asserted is (at least) true and regard a denial as correct just in case the sentence denied is (at least) false.  Officially, letting every sentence $A$ receive a value $v(A) \subseteq \{ 1, 0\}$, we may define the notion of correctness as follows:

\begin{quote}
\textbf{Correctness:} Asserting $A$ is \textit{correct}, relative to some valuation $v$, just in case $1 \in v(A) $.  Denying $A$ is \textit{correct}, relative to $v$, just in case $0 \in v(A)$.
\end{quote} 

\noindent Referring to assertions or denials generally as linguistic \textit{moves} one might make, we now define bilateral validity as follows:

\begin{quote}
 \textbf{Bilateral Validity:} An argument of the form $\Gamma \vdash \Delta$ is \textit{bilaterally valid}, relative to a set of valuations $V$, $\Gamma \vDash_{B_V} \Delta$, just in case there is no $v\in V$ such that all of the moves in $\Gamma$ are correct and all of the moves in $\Delta$ are incorrect.
 \end{quote}

\noindent This is just the notion of validity appealed to by Smiley and Rumfitt, which Rumfitt calls ``Smiley Consequence,'' generalized to multiple conclusions and 4-valued valuations.\footnote{See \cite{Blasio2017} for a systematic treatment of this notion of consequence under the heading of ``Two-Dimensional Consequence.''}

Now, just as we formulated the proof rules for these connectives schematically, we can also schematically state the Dunn semantic clauses relative to which these rules are sound and complete.    To do this, let us define the \textit{correctness function} as follows: 

  \singlespacing
\begin{quote}
\textbf{Correctness Function:} The \textit{correctness function} $[\cdot]$ is a function from $\{+, -\}$ to $\{1, 0\}$ mapping $+$ to $1$ and $-$ to $0$.  
\end{quote}

\noindent   The expression $[\textit{\textbf{a}}]$ can be read as ``the value that would make stance \textit{\textbf{a}} correct.''  This lets us define correctness, relative to a valuation $v$, abstractly as follows:
 \begin{quote}
\textbf{Correctness:} Taking some stance $\textit{\textbf{a}}$ towards some sentence $A$, $\textit{\textbf{a}} \langle A \rangle$, is \textit{correct}, relative to some valuation $v$, just in case $[\textit{\textbf{a}}] \in v(A) $.  
\end{quote}
 
 \noindent Having defined this notation, the semantic clauses for all of the binary connectives can be stated in a single schematic clause as follows:
 
\begin{quote}
\textbf{Aggregative Connectives:}

\begin{quote}
$ v(A \circ B) \ni
    \begin{cases}
      [\textit{\textbf{c}}], & \text{iff}\ [\textit{\textbf{a}}] \in v(A)   \text{ and } [\textit{\textbf{b}}]  \in v(B) \\
      [\textit{\textbf{c}}^*], & \text{iff}\ [\textit{\textbf{a}}^*]  \in v(A)  \text{ or } [\textit{\textbf{b}}^*] \in v(B)  
    \end{cases}$
\end{quote}

\end{quote}

\noindent One can check that this does indeed give the standard Dunn semantic clauses of familiar connectives, and in \citep{Simonelli2026} I schematically establish soundness and completeness relative to the proof rules introduced in Section 2. 

For the constructive connectives, let $\langle W,\leq\rangle$ be a Kripke frame, where $\leq$ is a preorder on $W$. A Dunn-Kripke valuation assigns each atomic formula $p$ a value $v_w(p)\subseteq\{1,0\}$ at every $w\in W$, subject to persistence: if $w\leq w'$, then $v_w(p)\subseteq v_{w'}(p)$. We take the above clause for the aggregative connectives to contain a world index, indicating that this semantic clause holds at every point.  For the Nelson-style connectives, we can formulate a schematic clause as follows: 

\begin{quote}
\textbf{Nelson-Style Connectives}

\begin{quote}
$v_w(A\triangleright B)\ni
\begin{cases}
[\signc],&
\text{iff }[\signa]\in v_w(A)
\text{ and }[\signb]\in v_w(B),
\\[1mm]
[\signc^*],&
\text{iff for every }w'\geq w,\ \text{ if } [\signa]\in v_{w'}(A)
\text{ then }
[\signb^*]\in v_{w'}(B).
\end{cases}$
\end{quote}

\end{quote}

\noindent Finally, Strong Negation is De Morgan Negation, given by the following Dunn semantic clause:

\begin{quote}
\textbf{Strong Negation:} For either $\signc\in\{+,-\}$, $[\signc]\in v_w(\csim A)$ iff $[\signc^*]\in v_w(A)$.
\end{quote}

\noindent It is straightforward to verify that the schematic proof rules given above are sound for these semantic clauses.\footnote{It is likewise straightforward to schematically establish completeness, but soundness is all that is needed for the following argument.}

\section{Defining Strong Negation with Two Connectives}

Section 3 established that a single constructive Nelson-style Sheffer stroke defines strong negation and thus the constructive conditional as well as the other members of the Nelson family.  The same can be shown for the Nelson-style Peirce's arrow (indeed, the very same formula defines negation).  Still, though we can define strong negation and the other Nelson-style connectives with such connectives, we are not able to define the full language of N4.  To obtain the full language of N4, we need an aggregative connective, from which conjunction and disjunction can be recovered.\footnote{I omit the proof for reasons of space, but the claim can be established by a preservation argument on the three-point fork frame \(r\leq u,r\leq v\), with \(u\) and \(v\) incomparable. The collection $\{\varnothing,\{u\},\{v\},\{r,u,v\}\}$ of upward-closed sets is closed under intersection and Heyting implication, and hence is preserved by every Nelson-style connective, but it is not closed under union, since \(\{u\}\cup\{v\}=\{u,v\}\). Suitable valuations therefore show that neither aggregative conjunction nor aggregative disjunction is definable from the Nelson-style connectives, even in the presence of strong negation. Since strong negation makes all aggregative connectives interdefinable, no aggregative connective is definable from the Nelson-style family.}  Now, once strong negation is available, every member of each family defines every other member of that family: applying strong negation to either argument changes the corresponding input sign, while applying strong negation to the entire compound changes its output sign.  So, the Nelson-style Sheffer stroke or Peirce's arrow, along with any of the aggregative connectives suffices.  Likewise the aggregative Sheffer stroke or Peirce's arrow, along with any of the Nelson-style connectives suffices.  The question remains of whether there are any \textit{other} selections of a single aggregative connective and a single Nelson-style connective that suffice to define strong negation.   

To answer the above question, let us first introduce a bit of notation.  To refer to different instances of the above schemas, while still potentially remaining at the schematic level, let us use a superscript $x,y,z$ on a connective symbol to indicate that the signs $x,y,z$ are substituted for $\signa,\signb,\signc$, respectively, in its rule schema and corresponding semantic clause. Each superscript entry may be a sign or a variable ranging over signs. Thus, while we might refer to the Nelson conditional as $\triangleright^{+, -, -}$, we might also write $\triangleright_j=\triangleright^{\signd,\signe,\signf}$ to schematically refer to a Nelson-style connective where $\signd$, $\signe$, $\signf$ are variables over $\{+,-\}$ corresponding to $\signa$, $\signb$, and $\signc$ in the above Nelson schema. We can now provide an exhaustive characterization of such two-connective pairs:

\begin{quote}
\textbf{Proposition 2:} Let $\circ_i=\circ^{\signa,\signb,\signc}$ be an aggregative connective and let $\triangleright_j=\triangleright^{\signd,\signe,\signf}$ be a Nelson-style connective.  Then, $\{\circ_i, \triangleright_j\}$ defines strong negation just in case at least one of the following conditions holds:
\begin{enumerate}[label=(\roman*)]
\item $\signa = \signb = \signc^*$
\item $\signd = \signe = \signf^*$
\item $\signa \neq \signb, \signd \neq \signe$ and $\signc = \signf^*$ 
\end{enumerate}
\end{quote}

\textit{Proof:} Let us first establish the sufficiency of the conditions purely proof-theoretically.  We have already established (i) and (ii), since these are just the cases in which one of the connectives is the (aggregative or Nelson-style) Sheffer stroke or Peirce's arrow.  In fact, though we just showed the proofs for the Sheffer stroke, they can be schematized so that they also show that negation can be defined by the same formula for Peirce's arrow.  For (iii), we can define strong negation as follows:

\begin{quote}
$\csim A :=
\begin{cases}
(A \triangleright_j A) \circ_i A,&\text{if }\textbf{\textit{a}}=\textit{\textbf{c}},\\
A \circ_i  (A \triangleright_j A)&\text{if } \textbf{\textit{b}} = \textbf{\textit{c}}
\end{cases}
$
\end{quote}

\noindent  To see how this definition works, consider first the diagonal formula $A\triangleright_j A$.  Since \(\signd\neq\signe\), we have \(\signe^*=\signd\), and so $\signf^*\langle A\triangleright_j A\rangle$ is always automatically derivable by Reflexivity, and since \(\signc=\signf^*\), this means that  $\signc \langle A\triangleright_j A\rangle$ is always derivable given Reflexivity:

\[
\infer[{\signf^*}_{\triangleright_j\mathrm R (\signc = \signf^*)}]
      {\Gamma\vdash
       \signc\langle A\triangleright_j A\rangle}
      {
       \infer[_{\text{\rm C. Reflex.} (\signe^*=\signd)}]
             {\Gamma,\signd\langle A\rangle
              \vdash\signe^*\langle A\rangle}
             {}
      }.
\]

\noindent By contrast, \(\signc^*=\signf\), so the \(\signc^* \langle A\triangleright_j A \rangle\) requires both \(\signd\langle A\rangle\) and \(\signe\langle A\rangle\). Since \(\signd\neq\signe\), one of these signs is $\signc$.  So, we have: 
\[
\infer[{\signf}_{\triangleright_j\mathrm L (\signc^* = \signf)}]
      {\signc^*\langle A\triangleright_j A\rangle
       \vdash\signc\langle A\rangle}
      {
       \infer[_{\text{\rm C. Reflex.} (\signd=\signc \text{ or } \signe = \signc)}]
             {\signd\langle A\rangle,
              \signe\langle A\rangle
              \vdash\signc\langle A\rangle}
             {}
      }.
\]

\noindent So, $A\triangleright_j A$ is such that one signed instance, $\signc \langle A\triangleright_j A \rangle$, is always derivable, and the other signed instance $\signc^* \langle A\triangleright_j A \rangle$, derives both $\signc\langle A\rangle$ and $\signc^* \langle A\rangle$.   Having used our Nelson-style connective to form a sentence with this inferential behavior, our aggregative connective, with the condition $\signa \neq \signb$, enables us to derive the four toggle sequents.  

We will just consider one case, as the other is analogous.  So, suppose that \(\signb=\signc\), and hence, since  \(\signa\neq\signb\), that \(\signa=\signc^*\). In this case, we define $\csim A:= A \circ_i  (A \triangleright_j A)$ and derive the four toggle sequents as follows:

\begin{center}
	\begin{minipage}{\linewidth}
		\centering

\begin{multicols}{2}

$\infer[{\signc}_{\circ_i\mathrm R (\signa \neq \signb)}]
      {\signc^*\langle A\rangle
       \vdash
       \signc\langle
       A\circ_i(A\triangleright_j A)
       \rangle}
      {
       \infer[_\text{\rm C. Reflex.}]
             {\signc^*\langle A\rangle
              \vdash\signc^*\langle A\rangle}
             {}
       &
       \infer[{\signf^*}_{\triangleright_j\mathrm R (\signc = \signf^*)}]
             {\signc^*\langle A\rangle
              \vdash
              \signc\langle A\triangleright_j A\rangle}
             {
              \infer[_\text{\rm C. Reflex.}]
                    {\signc^*\langle A\rangle,
                     \signd\langle A\rangle
                     \vdash\signd\langle A\rangle}
                    {}
             }
      }$

\columnbreak 

$\infer[{\signc}_{\circ_i\mathrm L (\signa \neq \signb)}]
      {\signc\langle
       A\circ_i(A\triangleright_j A)
       \rangle
       \vdash\signc^*\langle A\rangle}
      {
       \infer[_\text{\rm C. Reflex.}]
             {\signc^*\langle A\rangle,
              \signc\langle A\triangleright_j A\rangle
              \vdash\signc^*\langle A\rangle}
             {}
      }$

\end{multicols}
\begin{multicols}{2}

$\infer[{\signc^*}_{\circ_i\mathrm R (\signa \neq \signb)}]
      {\signc\langle A\rangle
       \vdash
       \signc^*\langle
       A\circ_i(A\triangleright_j A)
       \rangle}
      {
       \infer[_\text{\rm C. Reflex.}]
             {\signc\langle A\rangle
              \vdash
              \signc\langle A\rangle,
              \signc^*\langle A\triangleright_j A\rangle}
             {}
      }$

\columnbreak 

$\infer[{\signc^*}_{\circ_i\mathrm L (\signa \neq \signb)}]
      {\signc^*\langle
       A\circ_i(A\triangleright_j A)
       \rangle
       \vdash\signc\langle A\rangle}
      {
       \infer[_\text{\rm C. Reflex.}]
             {\signc\langle A\rangle
              \vdash\signc\langle A\rangle}
             {}
       &
       \infer[{\signf}_{\triangleright_j\mathrm L (\signc^* = \signf)}]
             {\signc^*\langle A\triangleright_j A\rangle
              \vdash\signc\langle A\rangle}
             {
              \infer[_{\text{\rm C. Reflex.}  (\signd=\signc \text{ or } \signe = \signc)}]
                    {\signd\langle A\rangle,
                     \signe\langle A\rangle
                     \vdash\signc\langle A\rangle}
                    {}
             }
      }$
\end{multicols}
\end{minipage}
\end{center}

\noindent Given Cut, these interderivabilities yield all four rules for strong negation.  If instead \(\signa=\signc\), then \(\signb=\signc^*\), and we define $\csim A:=(A\triangleright_j A)\circ_i A$ and the four derivations are exactly symmetric: the positions of the premises are merely exchanged, corresponding to the exchange in the arguments of the formula.  This establishes the sufficiency of the condition. 

We have established sufficiency proof-theoretically. For necessity, we use the semantics together with soundness.  Let us start with a basic observation connecting what we'll call ``common diagonal fixed values'' with proof-theoretic non-definability.  Consider the one-point Kripke frame \(W=\{w\}\), with \(w\leq w\). Call a classical value $ x\in\{\{1\},\{0\}\}$ a \emph{diagonal fixed value} of a binary connective \(C\) just in case $C(x,x)=x.$ Suppose that \(\circ_i\) and \(\triangleright_j\) have a common diagonal fixed value \(x=\{[\signg]\}\). Assign an atomic formula \(p\) the value $v_w(p)=\{[\signg]\}.$ A straightforward induction on formula complexity then shows that every term \(T(p)\) constructed from \(p\) using only \(\circ_i\) and \(\triangleright_j\) has the same value: $v_w(T(p))=\{[\signg]\}.$ Now, if such a term \(T(p)\) defined strong negation proof-theoretically, then, in particular, the following toggle sequent $\signg\langle T(p)\rangle \vdash \signg^*\langle p\rangle$ would be derivable.  Under the valuation just described, however, \(\signg\langle T(p)\rangle\) is correct, whereas \(\signg^*\langle p\rangle\) is incorrect. The sequent is therefore bilaterally invalid and, by soundness, is not derivable. Consequently, no term constructed from \(\circ_i\) and \(\triangleright_j\) can define strong negation whenever the two connectives have a common diagonal fixed value. We will now show that if none of conditions (i)-(iii) are met,  \(\circ_i\) and \(\triangleright_j\) must have such a common diagonal fixed value.  

A direct inspection of the semantic clause for aggregative connectives yields the following classification of the diagonal fixed-value behavior of an aggregative connective $\circ^{\signa,\signb,\signc}$:

\begin{enumerate}
\item if $\signa = \signb = \signc^*$, it fixes neither
classical signed value;
\item if $\signa = \signb = \signc$, it fixes both classical
signed values;
\item if $\signa \neq \signb$, it fixes exactly $\{[\signc^*]\}$
\end{enumerate}

\noindent Since we are considering a one-point frame and a classical valuation, exactly the same classification applies to a Nelson-style connective  $\triangleright^{\signd,\signe,\signf}$, with $\signd, \signe, \signf$ in place of $\signa, \signb, \signc$.    To see why, note that, on the one-point frame, its $\signf^*$-condition reduces to $[\signd]\in v(A)\Rightarrow [\signe^*]\in v(B)$. Moreover, on either classical value, exactly one of $[\signd]$ and $[\signd^*]$ is present. Hence, $[\signd]\in v(A)\Rightarrow [\signe^*]\in v(B)$ is equivalent, on these values, to $[\signd^*]\in v(A)$ or $[\signe^*]\in v(B)$. This has exactly the same form as the complementary condition for an aggregative connective. Consequently:

\begin{enumerate}
\item if $\signd=\signe=\signf^*$, it fixes neither classical value;
\item if $\signd=\signe=\signf$, it fixes both classical values;
\item if $\signd\neq\signe$, it fixes exactly $\{[\signf^*]\}$.
\end{enumerate}

\noindent Now, suppose that none of (i)--(iii) holds.  Once again, these are the following conditions: 

\begin{enumerate}[label=(\roman*)]
\item $\signa = \signb = \signc^*$
\item $\signd = \signe = \signf^*$
\item $\signa \neq \signb, \signd \neq \signe$ and $\signc = \signf^*$ 
\end{enumerate}

\noindent Since (i) fails, \(\circ_i\) is not of the first kind; since (ii) fails, \(\triangleright_j\) is not of the first kind.  If either connective is of the second kind, it fixes both classical signed values. The other connective, not being of the first kind, fixes at least one of them. Hence the two connectives have a common fixed classical signed value.  The only remaining possibility is that $\signa \neq \signb$ and $\signd \neq \signe$ . The unique fixed values of the two connectives are then, respectively, $\{[\signc^*]\}$  and $\{[\signf^*]\}$.  Since (iii) fails, $\signc \neq\signf^*$, and hence,  $\signc^* = \signf^*$, and thus $\{[\signc^*]\} = \{[\signf^*]\}$.  So the two connectives again have a common fixed value $\{[\textit{\textbf{g}}]\}$.   By the fixed-value observation above, no term constructed from \(\circ_i\) and \(\triangleright_j\) can define strong negation proof-theoretically. This establishes necessity.

Thus, \(\{\circ_i,\triangleright_j\}\) defines strong negation just in case at least one of (i)--(iii) holds. $\Box$

The above proof establishes precisely the conditions an aggregative/Nelson-style pair must meet if it is to be able to define strong negation.  We have already considered the pairs of connectives that meet (i) or (ii).  Let us, now, consider the pairs of connectives that do not meet (i) or (ii), but which meet (iii).  There are eight such pairs:

\begin{table}[H]
\centering
\renewcommand{\arraystretch}{1.4}
\begin{tabular}{c c l}
\hline
\textbf{Aggregative connective}
&
\textbf{Nelson-style connective}
&
\multicolumn{1}{c}{\textbf{Definition of \(\csim A\)}}
\\
\hline

\(\rightarrow_\circ\)
&
\(\starimp_\triangleright\)
&
\(\csim A:=A\rightarrow_\circ
  (A\starimp_\triangleright A)\)
\\

\(\rightarrow_\circ\)
&
\(\coimp_\triangleright\)
&
\(\csim A:=A\rightarrow_\circ
  (A\coimp_\triangleright A)\)
\\

\(\starcoimp_\circ\)
&
\(\starimp_\triangleright\)
&
\(\csim A:=
  (A\starimp_\triangleright A)
  \starcoimp_\circ A\)
\\

\(\starcoimp_\circ\)
&
\(\coimp_\triangleright\)
&
\(\csim A:=
  (A\coimp_\triangleright A)
  \starcoimp_\circ A\)
\\

\(\starimp_\circ\)
&
\(\rightarrow_\triangleright\)
&
\(\csim A:=
  (A\rightarrow_\triangleright A)
  \starimp_\circ A\)
\\

\(\starimp_\circ\)
&
\(\starcoimp_\triangleright\)
&
\(\csim A:=
  (A\starcoimp_\triangleright A)
  \starimp_\circ A\)
\\

\(\coimp_\circ\)
&
\(\rightarrow_\triangleright\)
&
\(\csim A:=A\coimp_\circ
  (A\rightarrow_\triangleright A)\)
\\

\(\coimp_\circ\)
&
\(\starcoimp_\triangleright\)
&
\(\csim A:=A\coimp_\circ
  (A\starcoimp_\triangleright A)\)
\\

\hline
\end{tabular}
\end{table}

\noindent Of course, most of these pairs contain at least one unusual connective, not included in the standard connective packages of standard logics.  Note, however, that in this bilateral context, there is no in principle reason to regard these unusual connectives as second-class citizens.  For instance, the constructive anti-conditional $\starimp_\triangleright$ is such that $A \starimp_\triangleright B$  is to be denied precisely when the assertion of $B$ can be concluded from the assertion of $A$.  In a bilateral context, where denial is treated as primitive, there is no reason to think that such a conditional is essentially defined as the negation of the constructive conditional $\rightarrow_\triangleright$.  Thus, it is perfectly reasonable to treat the anti-conditional as primitive and use the package shown in the first line to define both strong negation and the constructive conditional.

Though all of the conditionals displayed in the table above are reasonable candidates as primitive connectives in a bilateral context, the second line is still worth paying particular attention to, for both of the connectives belonging to this line \textit{are} standard: they are simply the material conditional and constructive co-implication.   Here are the rules for both connectives (here using $\supset$ for the material conditional and simply $\coimp$ with no subscript for the constructive co-implication):

\begin{center}

\begin{multicols}{2}

$
\infer[-_{\supset\mathrm L}]
      {\Gamma,\deny{A\supset B}\vdash\Delta}
      {\Gamma,\assert A,\deny B\vdash\Delta}
$

\columnbreak

$
\infer[-_{\supset\mathrm R}]
      {\Gamma\vdash\deny{A\supset B},\Delta}
      {\Gamma\vdash\assert A,\Delta
       &
       \Gamma\vdash\deny B,\Delta}
$

\end{multicols}

\begin{multicols}{2}

$
\infer[+_{\supset\mathrm L}]
      {\Gamma,\assert{A\supset B}\vdash\Delta}
      {\Gamma,\deny A\vdash\Delta
       &
       \Gamma,\assert B\vdash\Delta}
$

\columnbreak

$
\infer[+_{\supset\mathrm R}]
      {\Gamma\vdash\assert{A\supset B},\Delta}
      {\Gamma\vdash\deny A,\assert B,\Delta}
$

\end{multicols}

\end{center}

\begin{center}

\begin{multicols}{2}

$
\infer[+_{\coimprule\mathrm L}]
      {\Gamma,\assert{A\coimp B}\vdash\Delta}
      {\Gamma,\deny A,\assert B\vdash\Delta}
$

\columnbreak

$
\infer[+_{\coimprule\mathrm R}]
      {\Gamma\vdash\assert{A\coimp B},\Delta}
      {\Gamma\vdash\deny A,\Delta
       &
       \Gamma\vdash\assert B,\Delta}
$

\end{multicols}

\begin{multicols}{2}

$
\infer[-_{\coimprule\mathrm L}]
      {\Gamma,\deny{A\coimp B}\vdash\Delta}
      {\Gamma,\deny{A\coimp B}\vdash\deny A
       &
       \Gamma,\deny B\vdash\Delta}
$

\columnbreak

$
\infer[-_{\coimprule\mathrm R}]
      {\Gamma\vdash\deny{A\coimp B},\Delta}
      {\Gamma,\deny A\vdash\deny B}
$

\end{multicols}

\end{center}

\noindent It is perhaps surprising that this is in fact an expressively complete basis for N4, but the following proofs, which are an instance of the schematic proofs shown in the proof above, show that $A\supset (A\coimp A)$ does indeed define strong negation:

\small
\begin{center}
\begin{multicols}{2}

$
\infer[+_{\supset\mathrm R}]
{
  -\langle A\rangle
  \vdash
  +\langle
  A\supset
  (A\coimp A)
  \rangle
}
{
  \infer[_\text{C. Reflex.}]
  {
    -\langle A\rangle
    \vdash
    -\langle A\rangle,
    +\langle A\coimp A\rangle
  }
  {}
}
$

\columnbreak 

$
\infer[+_{\supset\mathrm L}]
{
  +\langle
  A\supset
  (A\coimp A)
  \rangle
  \vdash
  -\langle A\rangle
}
{
  \infer[_\text{C. Reflex.}]
  {
    -\langle A\rangle
    \vdash
    -\langle A\rangle
  }
  {}
  &
  \infer[+_{{\coimprule}_ L}]
  {
    +\langle A\coimp A\rangle
    \vdash
    -\langle A\rangle
  }
  {
    \infer[_\text{C. Reflex.}]
    {
      -\langle A\rangle,
      +\langle A\rangle
      \vdash
      -\langle A\rangle
    }
    {}
  }
}
$
\end{multicols}
\begin{multicols}{2}

$
\infer[-_{\supset\mathrm R}]
{
  +\langle A\rangle
  \vdash
  -\langle
  A\supset
  (A\coimp A)
  \rangle
}
{
  \infer[_\text{C. Reflex.}]
  {
    +\langle A\rangle
    \vdash
    +\langle A\rangle
  }
  {}
  &
  \infer[-_{\coimprule\mathrm R}]
  {
    +\langle A\rangle
    \vdash
    -\langle A\coimp A\rangle
  }
  {
    \infer[_\text{C. Reflex.}]
    {
      +\langle A\rangle,
      -\langle A\rangle
      \vdash
      -\langle A\rangle
    }
    {}
  }
}
$

\columnbreak 

$
\infer[-_{\supset\mathrm L}]
{
  -\langle
  A\supset
  (A\coimp A)
  \rangle
  \vdash
  +\langle A\rangle
}
{
  \infer[_\text{C. Reflex.}]
  {
    +\langle A\rangle,
    -\langle A\coimp A\rangle
    \vdash
    +\langle A\rangle
  }
  {}
}
$
\end{multicols}
\end{center}
\normalsize

\noindent Having defined strong negation, all of the other connectives of both families can be defined, given $\supset$ and $\coimp$.  

This classification suggests that the definability of strong negation is considerably less exceptional than Oddsson’s four-connective definition might initially suggest. Oddsson’s definition employs the full binary vocabulary typical of 2Int: two aggregative and two Nelson-style connectives. Once all the connectives generated by the two rule schemas are treated as equally legitimate, however, a majority of the possible two-connective packages suffice. Of the 64 packages containing one aggregative and one Nelson-style connective, 36 define strong negation.  28 do so because one member already defines negation individually, while 8 further pairs---despite neither member defining it individually---define it jointly.

\section{Wansing-Style Connexive Connectives}

Thus far, I've given two positive results concerning the definability of strong negation with Nelson-style constructive connectives.   Given the negative results of \cite{WansingNikiDrobyshevich2025} it is natural to expect that strong negation would be harder to come by if we consider \textit{connexive} conditionals of the sort proposed by Wansing \citeyearpar{Wansing2005} \citeyearpar{Wansing2016b}.  Surprisingly, we have striking positive definability results here as well, indeed, a direct analogue of Proposition 1.

To briefly provide a bit of background for these connectives, connexive logic is a contra-classical logic, characterized by the validation of certain non-theorems of classical logic, most notably, \textit{Aristotle's Theses} and \textit{Boethius’s Theses}.  In a bilateral system, we might put these as follows:

\begin{quote}
\textbf{Aristotle's Theses:} 
\begin{enumerate}
\item $\vdash - \langle A \rightarrow \csim A \rangle$ 
\item $\vdash - \langle \csim A \rightarrow A \rangle$ 
\end{enumerate}  

\textbf{Boethius's Theses:} 
\begin{enumerate}
\item $\vdash + \langle (A \rightarrow B) \rightarrow \csim (A \rightarrow \csim B) \rangle$
\item  $\vdash + \langle (A \rightarrow \csim B) \rightarrow \csim (A \rightarrow B) \rangle$
\end{enumerate}
\end{quote}

\noindent So, for any sentences $A$ and $B$, we're always committed to denying that $A$ implies its own negation (and that the negation of $A$ implies $A$), and we're always committed to asserting that if $A$ implies $B$, then it's not the case that $A$ implies the negation of $B$ (and if $A$ implies the negation of $B$, then it's not the case that $A$ implies $B$).  Wansing shows that, in a bilateral context, there is a very natural way of meeting the desiderata: simply take the rules for a Nelson-style conditional and tweak the denial condition so that it inferentially internalizes the transition from the assertion of the antecedent to the denial of the consequent.  Let us refer to the following as the ``Wansing conditional'' (the multiple conclusion rules, once again, are owed to \cite{Drobyshevich2026}):

\begin{center}

   \begin{multicols}{2}
   
    $
  \infer[+_{\rightarrow_\text{L}}]
  {\Gamma, \assert{A \rightarrow B} \vdash \Delta}
  {\Gamma, \assert{A \rightarrow B}
       \vdash \assert A
   &
   \Gamma, \assert B \vdash \Delta}
  $
  
    \columnbreak
  
    $
  \infer[+_{\rightarrow_\text{R}}]
  {\Gamma \vdash
       \assert{A \rightarrow B}, \Delta}
  {\Gamma, \assert A \vdash \assert B}
  $
  
   \end{multicols}

   \begin{multicols}{2}
   
    $
  \infer[-_{\rightarrow_\text{L}}]
  {\Gamma, \deny{A \rightarrow B} \vdash \Delta}
  {\Gamma, \deny{A \rightarrow B}
       \vdash \assert A
   &
   \Gamma, \deny B \vdash \Delta}
  $
  
    \columnbreak
  
    $
  \infer[-_{\rightarrow_\text{R}}]
  {\Gamma \vdash
       \deny{A \rightarrow B}, \Delta}
  {\Gamma, \assert A \vdash \deny B}
  $
  
   \end{multicols}

\end{center}

\noindent As you can see, the assertion condition for $A \rightarrow B$ is the same as that of the Nelson conditional: it internalizes the positive inferential relation between $A$ and $B$---that assertion of $B$ follows from the assertion of $A$.  Unlike the Nelson conditional, however, the denial condition for $A \rightarrow B$ internalizes the \textit{negative} inferential relation between $A$ and $B$---that the denial of $B$ follows from the assertion of $A$.  With these rules and the rules for strong negation, Aristotle's Theses and Boethius’s Theses are straightforwardly derivable.

Now, following the same generalized approach as before, let us consider the class of connexive Wansing-style connectives:

\begin{quote}
\textbf{Wansing-Style Connectives:}

\begin{center}

   \begin{multicols}{2}
   
    $
  \infer[{{\textbf{\textit{c}}}_{\blacktriangleright_\text{L}}}]{\Gamma,  \textbf{\textit{c}}  \langle A \blacktriangleright B \rangle \vdash \Delta}{\Gamma,   \textbf{\textit{c}}  \langle A \blacktriangleright B \rangle   \vdash \textbf{\textit{a}}  \langle  A \rangle  & \Gamma, \textbf{\textit{b}} \langle B \rangle \vdash \Delta}
  $
  
    \columnbreak
  
       $
  \infer[{{\textbf{\textit{c}}}_{\blacktriangleright_\text{R}}}]{\Gamma \vdash \textbf{\textit{c}} \langle A \blacktriangleright B \rangle, \Delta}{\Gamma, \textbf{\textit{a}}  \langle A \rangle \vdash  \textbf{\textit{b}} \langle B \rangle}
  $

   \end{multicols}

   \begin{multicols}{2}
   
    $
  \infer[{{\textbf{\textit{c}}^*}_{\blacktriangleright_\text{L}}}]{\Gamma,  \textbf{\textit{c}}^*  \langle A \blacktriangleright B \rangle \vdash \Delta}{\Gamma, \textbf{\textit{c}}^*  \langle A \blacktriangleright B \rangle \vdash \textbf{\textit{a}}  \langle  A \rangle   & \Gamma, \textbf{\textit{b}}^* \langle B \rangle \vdash \Delta}
  $
  
    \columnbreak
  
       $
  \infer[{{\textbf{\textit{c}}^*}_{\blacktriangleright_\text{R}}}]{\Gamma \vdash \textbf{\textit{c}}^* \langle A \blacktriangleright B \rangle, \Delta}{\Gamma, \textbf{\textit{a}}  \langle A \rangle \vdash  \textbf{\textit{b}}^* \langle B \rangle}
  $

   \end{multicols}

 \end{center}
\end{quote}

\noindent Unlike the previous two classes of connectives, where there were \textit{eight} distinct instances, corresponding to each of the different assignments of signs to \signa, \signb, and \signc, here there are only four distinct instances, since the assignments $(\signa, \signb, \signc)$ and $(\signa, \signb^*, \signc^*)$ collapse.  I thus choose the following four distinct representative connectives:

\begin{center}
\begin{multicols}{2}

$\mid_\blacktriangleright$: $\textbf{\textit{a}}=+, \textbf{\textit{b}}=+, \textbf{\textit{c}}=- $\\

$\rightarrow_\blacktriangleright$: $\textbf{\textit{a}}=+, \textbf{\textit{b}}=-, \textbf{\textit{c}}=-$

\columnbreak

$\downarrow_\blacktriangleright$: $\textbf{\textit{a}}=-, \textbf{\textit{b}}=-, \textbf{\textit{c}}=+$\\

$\coimp_\blacktriangleright$: $\textbf{\textit{a}}=-, \textbf{\textit{b}}=+, \textbf{\textit{c}}=+$

\end{multicols}
\end{center}

\noindent As before, we can provide the Kripke-Dunn style semantic clauses for all of these connectives schematically:

\begin{quote}
\textbf{Constructive Wansing-Style Connectives}

\begin{quote}
$v_w(A\blacktriangleright B)\ni
\begin{cases}
[\signc],&
\text{iff for every }w'\geq w,\ \text{ if } [\signa]\in v_{w'}(A)
\text{ then }
[\signb]\in v_{w'}(B).
\\[1mm]
[\signc^*],&
\text{iff for every }w'\geq w,\ \text{ if } [\signa]\in v_{w'}(A)
\text{ then }
[\signb^*]\in v_{w'}(B).
\end{cases}$
\end{quote}

\end{quote}

\noindent As before, it is easy to check soundness, which is all we will need for the following results.

Once again, our focus will be on the Wansing-style Sheffer stroke, \(\mid_{\blacktriangleright}\), which is given by the following rules (to reduce clutter, the subscript \(\blacktriangleright\) is suppressed throughout this section):

\begin{center}
\begin{multicols}{2}

$\infer[-_{\stroke\mathrm L}]
      {\Gamma,\deny{A\stroke B}\vdash\Delta}
      {\Gamma,\deny{A\stroke B}\vdash\assert A
       &
       \Gamma,\assert B\vdash\Delta}$

\columnbreak

$\infer[-_{\stroke\mathrm R}]
      {\Gamma\vdash\deny{A\stroke B},\Delta}
      {\Gamma,\assert A\vdash\assert B}$

\end{multicols}

\begin{multicols}{2}

$\infer[+_{\stroke\mathrm L}]
      {\Gamma,\assert{A\stroke B}\vdash\Delta}
      {\Gamma,\assert{A\stroke B}\vdash\assert A
       &
       \Gamma,\deny B\vdash\Delta}$

\columnbreak

$\infer[+_{\stroke\mathrm R}]
      {\Gamma\vdash\assert{A\stroke B},\Delta}
      {\Gamma,\assert A\vdash\deny B}$

\end{multicols}
\end{center}

\noindent We will now show the following:

\begin{quote}
\textbf{Proposition 3:} The system containing just the Wansing-style Sheffer stroke rules enables one to define strong negation, and thus all Wansing-style connectives.
\end{quote}

\textit{Proof:} Let us define strong negation by:

\begin{quote}
$\csim A :=
((A\stroke
      (A\stroke A)))
\stroke A$
\end{quote}

\noindent To show that this formula indeed defines \(\csim A\), we show that it does in fact yield the relevant toggle behavior of strong negation:

\begin{quote}
\small

$\deny A \vdash
\assert{
((A\stroke
      (A\stroke A)))
\stroke A}$:

\begin{center}
$\infer[+_{\stroke\mathrm R}]
{
  \deny A
  \vdash
  \assert{
  ((A\stroke
        (A\stroke A)))
  \stroke A}
}
{
  \infer[_\text{C. Reflex.}]
  {
    \deny A,
    \assert{
    A\stroke
    (A\stroke A)}
    \vdash
    \deny A
  }
  {}
}$
\end{center}

$\assert{
((A\stroke
      (A\stroke A)))
\stroke A}
\vdash \deny A$:

\begin{center}
$\infer[+_{\stroke\mathrm L}]
{
  \assert{
  ((A\stroke
        (A\stroke A)))
  \stroke A}
  \vdash
  \deny A
}
{
  \infer[+_{\stroke\mathrm R}]
  {
    \assert{
    ((A\stroke
          (A\stroke A)))
    \stroke A}
    \vdash
    \assert{
    A\stroke
    (A\stroke A)}
  }
  {
    \infer[-_{\stroke\mathrm R}]
    {
      \assert{
      ((A\stroke
            (A\stroke A)))
      \stroke A},
      \assert A
      \vdash
      \deny{A\stroke A}
    }
    {
      \infer[_\text{C. Reflex.}]
      {
        \assert{
        ((A\stroke
              (A\stroke A)))
        \stroke A},
        \assert A,
        \assert A
        \vdash
        \assert A
      }
      {}
    }
  }
  &
  \infer[_\text{C. Reflex.}]
  {
    \deny A\vdash\deny A
  }
  {}
}$
\end{center}

$\assert A\vdash
\deny{
((A\stroke
      (A\stroke A)))
\stroke A}$:

\begin{center}
$\infer[-_{\stroke\mathrm R}]
{
  \assert A
  \vdash
  \deny{
  ((A\stroke
        (A\stroke A)))
  \stroke A}
}
{
  \infer[_\text{C. Reflex.}]
  {
    \assert A,
    \assert{
    A\stroke
    (A\stroke A)}
    \vdash
    \assert A
  }
  {}
}$
\end{center}

$\deny{
((A\stroke
      (A\stroke A)))
\stroke A}
\vdash\assert A$:

\begin{center}
$\infer[-_{\stroke\mathrm L}]
{
  \deny{
  ((A\stroke
        (A\stroke A)))
  \stroke A}
  \vdash
  \assert A
}
{
  \infer[+_{\stroke\mathrm R}]
  {
    \deny{
    ((A\stroke
          (A\stroke A)))
    \stroke A}
    \vdash
    \assert{
    A\stroke
    (A\stroke A)}
  }
  {
    \infer[-_{\stroke\mathrm R}]
    {
      \deny{
      ((A\stroke
            (A\stroke A)))
      \stroke A},
      \assert A
      \vdash
      \deny{A\stroke A}
    }
    {
      \infer[_\text{C. Reflex.}]
      {
        \deny{
        ((A\stroke
              (A\stroke A)))
        \stroke A},
        \assert A,
        \assert A
        \vdash
        \assert A
      }
      {}
    }
  }
  &
  \infer[_\text{C. Reflex.}]
  {
    \assert A\vdash\assert A
  }
  {}
}$
\end{center}

\normalsize
\end{quote}

\noindent Given Cut, we can derive the standard positive and negative left and right rules for toggle negation.

Having defined strong negation, we may define the constructive connexive conditional as follows: 

\begin{quote}
$A\rightarrow_{\blacktriangleright}B := A\stroke\csim B$
\end{quote}

\noindent  The same explanation for why this works applies as before.  We can likewise define the other Wansing-style connectives. So, the Wansing-style Sheffer stroke defines all of the Wansing-style connectives considered here. $\Box$

I hope the generalized approach that I have adopted here has made it clear why it is reasonable to conceive of the connective $\mid_\blacktriangleright$, defined by the rules above, as a connexive Sheffer stroke.  However, whereas the Nelson-style Sheffer stroke is novel to this paper (at least, as far as I'm aware), what I am calling the ``Wansing-style Sheffer stroke'' has in fact recently been discussed by \cite{ShramkoWansing2025}, who refer to this connective as \textit{connexive exclusion}.\footnote{Somewhat confusingly, they also refer to this connective as a kind of ``co-implication,'' though based on a distinct notion of duality as compared to the typical connexive co-implication $\coimp$, which figures in the Bi-Connexive logic 2C.}   Now, technically, they consider the connective $\mathbin{\reflectbox{$\mid_{\blacktriangleright}$}}$, which is just the connective $\mid_\blacktriangleright$ with the argument places reversed, such that $A \mathbin{\reflectbox{$\mid_{\blacktriangleright}$}} B = B \mid_\blacktriangleright A$, but, once again, this is a minor difference in notational convention and nothing significant hangs on it.  More interestingly, the symbol they use to express connexive exclusion is $\nleftarrow$.  The difference in notation here is explicable from the fact that, whereas we picked $\mid_\blacktriangleright$, treating our representative connective as a kind of Sheffer stroke, we could have just as well picked as our representative connective $\starimp_\blacktriangleright$, since rules associated with the sign assignments for $\mid$ and $\starimp$ collapse in this connexive context.  The Sheffer stroke, in general, however, can be understood as a kind of incompatibility exclusion operator---expressing the exclusion between two sentences---and so it is indeed natural to think of ``connexive exclusion'' as a kind of Sheffer stroke.  What Shramko and Wansing don't observe is that, like the standard Sheffer stroke, it is an expressively complete basis for the strong negation and implication fragment of Wansing's connexive logic.  Likewise, the Wansing-style Peirce's arrow defines strong negation. In general, a Wansing-style connective $\blacktriangleright^{\signa, \signb, \signc}$ defines strong negation just in case $\signb \neq \signc$.

Now, we saw, in the case of Nelson-style connectives and the aggregative connectives, that there were many pairs such that, although neither connective defines strong negation by itself, the two define it together. Nothing analogous occurs when the Nelson-style connective is replaced by a Wansing-style connective:

\begin{quote}
\textbf{Proposition 4.}
Let
\(\circ_i=\circ^{\signa,\signb,\signc}\)
be an aggregative connective and let
\(\blacktriangleright_j
=\blacktriangleright^{\signd,\signe,\signf}\)
be a Wansing-style connective. If neither connective defines strong
negation by itself, then they do not define strong negation jointly.
Equivalently, if
\[
\signa=\signb=\signc^*
\]
fails and
\[
\signe=\signf,
\]
then
\(\{\circ_i,\blacktriangleright_j\}\) does not define strong negation.
\end{quote}

\textit{Proof:}
Consider the one-point Kripke frame \(W=\{w\}\). Since
\(\signa=\signb=\signc^*\) fails, either
\(\signa=\signb=\signc\) or
\(\signa\neq\signb\). Let
\[
\signg = 
\begin{cases}
\signc,
&\text{if }\signa=\signb=\signc,\\
\signc^*,
&\text{if }\signa\neq\signb.
\end{cases}
\]
Assign every atomic formula the value \(\{[\signg]\}\). We show by
induction on formulas that every formula constructed using
\(\circ_i\) and \(\blacktriangleright_j\) contains
\([\signg]\).  For the aggregative connective, if
\(\signa=\signb=\signc\), its
\(\signc=\signg\)-condition is preserved. If
\(\signa\neq\signb\), exactly one of
\(\signa^*\) and \(\signb^*\) is
\(\signc^*=\signg\), so its
\(\signc^*=\signg\)-condition is preserved.  The Wansing-style connective preserves either sign when
\(\signe=\signf\). If \(\signg=\signf\), its
\(\signf\)-condition is satisfied because its consequent has sign
\(\signe=\signf\). If
\(\signg=\signf^*\), its \(\signf^*\)-condition is satisfied
because its consequent has sign
\(\signe^*=\signf^*\).  Thus every formula in the language contains \([\signg]\). 

Now, let \(T(p)\) be any term constructed from \(\circ_i\) and \(\blacktriangleright_j\). If \(T(p)\) defined strong negation proof-theoretically, then the following toggle sequent would be derivable: $\signg\langle T(p)\rangle \vdash \signg^*\langle p\rangle.$  Under the valuation just described, \(\signg\langle T(p)\rangle\) is correct, since \([\signg]\in v_w(T(p))\), whereas \(\signg^*\langle p\rangle\) is incorrect, since \(v_w(p)=\{[\signg]\}\). The sequent is therefore bilaterally invalid and, by soundness, is not derivable. Hence no term constructed using \(\circ_i\) and \(\blacktriangleright_j\) defines strong negation. \(\Box\)

\section{Conclusion}

I have presented here three distinct bilateral Sheffer strokes: the aggregative stroke $\mid_\circ$, the Nelson-style stroke $\mid_\triangleright$, and the Wansing-style stroke $\mid_\blacktriangleright$.  Each stroke by itself suffices to define strong negation and thus all of the connectives belonging to the families to which they belong.  Conceptually, the three strokes may be understood as incompatibility operators, with the latter two obtained by two successive bilateral intensionalizations of the incompatibility relation expressed by the previous operator. The aggregative stroke expresses merely extensional incompatibility: it is to be asserted when at least one argument is denied and it is to be denied when both arguments are asserted. The Nelson-style stroke retains this aggregative denial condition but replaces the assertion condition with an intensional, constructive requirement: that denial of the consequent follow from assertion of the antecedent. The Wansing-style stroke retains this intensional assertion condition while also intensionalizing denial: it is to be denied when assertion of the consequent follows from assertion of the antecedent. The incompatibility relation codified by these strokes \textit{can} be expressed with a corresponding positive conditional and strong negation: in each case, $A \mid_x B$ is equivalent to $A \rightarrow_x \sim B$.  However, in a bilateral setting there is reason to regard the incompatibility operator as more fundamental than the corresponding positive conditional. For the intensional connectives, whereas the assertion of a positive conditional expresses the \textit{unilateral} relation between the \textit{assertion} of one sentence and the \textit{assertion} of another, the assertion of the incompatibility operator directly expresses the \textit{bilateral} relation between the \textit{assertion} of one sentence and the \textit{denial} of another.  In a bilateral context, the expression of such a relation need not be mediated by negation.  On the contrary, as we've seen, we can define the operator of negation in terms of such relations.

Often, the Sheffer stroke is treated as little more than a technical curiosity: a way of arriving at an expressively parsimonious system for a logic with little bearing on how that logic should actually be conceptualized.  Indeed, towards the end of Shramko and Wansing's paper on connexive exclusion, they discuss the Sheffer stroke as an ``example of a case speaking against expressive parsimony,'' writing:

\begin{quote}
Although the set $\{\mid\}$ is functionally complete for Boolean propositional logic, the logic is usually presented by using the set of connectives $\{\neg, \wedge, \vee, \supset\}.$ The choice of connectives matters, and it matters not only for practical reasons such as succinctness, readability, or processability of notation. The choice of connectives can be a source of inspiration for further (in particular, philosophical) questions and investigations [\ldots] \citeyearpar[3396]{ShramkoWansing2025}
\end{quote}

\noindent They make this remark in the context of defending the treatment of connexive exclusion as a primitive binary connective in a connexive system, to be placed alongside strong negation and connexive implication, even though it can be defined in terms of the two. The pleasing irony here, as we've now seen, is that connexive exclusion is reasonably conceived as a Sheffer stroke itself, and can be used to define both strong negation and connexive implication.  The expressive parsimony of the one-connective bilateral systems is not accidental, but rather a consequence of the fact that the one-connective presentation captures the bilateral core of the logic.  

The same can be said, albeit in a different way, about the generalized notation through which these results were reached. I have developed and deployed this notation in a number of papers on bilateralism over the past few years \citeyearpar{Simonelli2024,Simonelli2025,Simonelli2026}, and I have found working in this notation to be very productive, not just technically, but conceptually as well. It is easy to think that this generalized notation is a way of arriving at an expressively parsimonious metalanguage: a compact way of stating systems consisting in independent connective rules with little bearing on how those systems should be conceptualized.  I hope to have shown here---particularly with the two-connective results, both positive and negative---that this approach can be used to arrive at general results about bilateral systems, articulating structural relations between kinds of connective rules that would be obscured by considering each set of connective rules independently.  Once again, the expressive parsimony of the generalized notation is not accidental, but rather a consequence of the fact that it captures the common bilateral core of the connective rules formulated with the use of it. In this sense, the choice of metalanguage notation, no less than the choice of object-language connectives, can be a source of inspiration for further philosophical and logical questions and investigations. \\

\noindent \textbf{Acknowledgments:} Thanks to Sara Ayhan, Sergey Drobyshevich, Bob Brandom, Ulf Hlobil, and an anonymous referee of this journal for comments.  

OpenAI's ChatGPT (5.6 Sol) was used during the exploratory stages of this research (in July 2026) to help search the space of possible connectives and relations among them, find relevant literature, and assist in the preliminary construction of the proofs of Propositions 2 and 4 (as well as the proof omitted but mentioned in footnote 13). I subsequently worked through and independently verified all of the outputs on which I drew in writing this paper, wrote all of the text contained in this paper myself, and take full responsibility for any errors this paper may contain. After the acceptance of the paper (in November 2026), I had GPT-6 Astra double check the paper's mathematical arguments, at which point it found an error in one of the works originally cited (which I independently verified), now noted in footnote 9.

\singlespacing
\newif\ifprojectbibliography
\projectbibliographyfalse
\ifprojectbibliography
\bibliographystyle{plainnat}
\bibliography{bochumbib}
\else

\fi
\onehalfspacing

\end{document}